\documentclass[12pt]{article}

\usepackage{graphicx}
\usepackage{amsmath, amsthm}
\usepackage{amsfonts}
\usepackage{amssymb}
\usepackage{url}
\usepackage{verbatim}
\usepackage{tikz}
\usepackage{color}
\definecolor{red}{rgb}{1,0,0}

\definecolor{blu}{rgb}{0,0,0}  
\def\blu{\color{blu}}
\usepackage[mathlines]{lineno}
\usepackage{mathrsfs}
\definecolor{qqqqff}{rgb}{0.,0.,1.}

\def\noi{\noindent}

\newcommand{\nul}{\operatorname{null}}
\newcommand{\rank}{\operatorname{rank}}
\makeatletter
\newcommand{\eqnum}{\refstepcounter{equation}\textup{\tagform@{\theequation}}}
\makeatother

\newlength{\sdgwidth}
\def\sdiv#1{\settowidth{\sdgwidth}{$#1$}
\mathop{\ooalign{$\mkern1.6mu\overline{\protect\raisebox{0.4pt}{$
\phantom{\mkern-1.6mu#1\kern-0.5\sdgwidth\mkern-.8mu}$}}$\hidewidth\protect\cr$
\kern0.5\sdgwidth\mkern0.8mu\mathring{\protect\raisebox{-.25pt}{$
\phantom{\kern-0.5\sdgwidth#1\kern-0.5\sdgwidth}$}}
\mkern-0.8mu\kern0.5\sdgwidth$\protect\cr$#1$\protect\cr$
\hidewidth\overline{\protect\raisebox{0.4pt}{$
\phantom{\kern-0.5\sdgwidth\mkern-2.2mu#1}$}}$}}\nolimits}

\newcommand{\R}{\mathbb{R}}

\newcommand{\Z}{\operatorname{Z}}

\newcommand{\Ad}{\mathscr {A}}

\newcommand{\SG}{\mathcal{S}(G)}

\newcommand{\sym}{\mathcal{S}}
\newcommand{\mr}{\operatorname{mr}}

\newcommand{\M}{\operatorname{M}}

\newtheorem{thm}{Theorem}[section]
\newtheorem{cor}[thm]{Corollary}
\newtheorem{dfn}[thm]{Definition}
\newtheorem{lem}[thm]{Lemma}
\newtheorem{prop}[thm]{Proposition}

\newtheorem{obs}[thm]{Observation}

\theoremstyle{definition}
\newtheorem{rem}[thm]{Remark}

\theoremstyle{definition}
\newtheorem{defn}[thm]{Definition}

\theoremstyle{definition}

\newcommand{\bit}{\begin{itemize}}
\newcommand{\eit}{\end{itemize}}
\newcommand{\ben}{\begin{enumerate}}
\newcommand{\een}{\end{enumerate}}
\newcommand{\beq}{\begin{equation}}
\newcommand{\eeq}{\end{equation}}
\newcommand{\bea}{\begin{eqnarray*}}
\newcommand{\eea}{\end{eqnarray*}}
\newcommand{\bpf}{\begin{proof}}
\newcommand{\epf}{\end{proof}\ms}
\newcommand{\ms}{\medskip}

\title{Zero forcing in iterated line digraphs}

\author{ Daniela Ferrero\thanks{Department of Mathematics, Texas State University, San Marcos, TX
    78666, U.S.A., \texttt{dferrero@txstate.edu}}\and Thomas Kalinowski\thanks{School of Science and Technology, University of New England, Armidale, NSW 2351, Australia, \texttt{tkalinow@une.edu.au}}\and Sudeep Stephen\thanks{School of Mathematical
    and Physical Sciences, The University of
    Newcastle, Callaghan, NSW 2308, Australia,  \texttt{sudeep.stephen@uon.edu.au}}}

\begin{document}

\date{}
\maketitle

\begin{abstract} 
Zero forcing is a propagation process on a graph, or digraph, defined in linear algebra to provide a bound for the minimum rank problem. Independently, zero forcing was introduced in physics, computer science and network science, areas where line digraphs are frequently used as models. Zero forcing is also related to power domination, a propagation process that models the monitoring of electrical power networks.

In this paper we study zero forcing in iterated line digraphs and provide a relationship between zero forcing and power domination in line digraphs. In particular, for regular iterated line digraphs we determine the minimum rank/maximum nullity, zero forcing number and power domination number, and provide constructions to attain them.  We conclude that regular iterated line digraphs present optimal minimum rank/maximum nullity, zero forcing number and power domination number, and apply our results to determine those parameters on some families of digraphs often used in applications.
\end{abstract}

\noi {\bf Keywords:} Minimum rank; Zero forcing; Power domination; Iterated line digraphs.  
\smallskip

\noi{\bf AMS subject classification} 05C20, 05C50, 05C76


\section{Introduction}
Zero forcing is a propagation process on a graph, independently introduced in linear algebra, physics, computer science and network science. In linear algebra, zero forcing was introduced in \cite{AIM08} to express a bound for the minimum rank problem, which consists of minimizing the rank of a matrix whose pattern of non-zero entries is determined by a given graph. The minimum rank problem appears frequently in engineering, where the order or the complexity of models can often be expressed as the rank of a matrix {\blu (see \cite{FHB04})}. In physics, zero forcing was introduced to study controllability of quantum systems (see \cite{BG07}); in computer science, it appears as the {\em fast-mixed search} model {\blu (see \cite{B13})} for some {\em pursuit-evasion} games {\blu (see \cite{P76})}; in network science, it models the spread of a disease over a population, or of an opinion in a social network {\blu (see \cite{DR})}. 

In addition to its intrinsic relation to minimum rank, zero forcing is closely related to power domination, a graph theory concept introduced in \cite{HHHH02} to optimize the monitoring process of electrical power networks. From the definitions of power domination and zero forcing, it follows that the closed out neighborhood of a power dominating set is a zero forcing set, and a stronger relationship between zero forcing and power domination was established in \cite{REUF2015}.

{\blu Zero forcing, minimum rank and power domination} are all {\it NP}-hard problems, as proven in {\blu \cite{A08}, \cite{BFS99},  and \cite{HHHH02}}, respectively. Thus, it is important to obtain bounds for the minimum rank, the zero forcing and the power domination numbers, as well as closed formulas to calculate them for families of graphs. Although {\blu power domination and zero forcing} where introduced on undirected graphs, they were extended to digraphs in {\blu \cite{A08} and   \cite{B09}}, respectively. Further results on zero forcing on digraphs can be found in {\blu \cite{B09},  \cite{FH11}, \cite{H10}, \cite{MZC14}  and \cite{TD15}.} Power domination in digraphs has not been so thoroughly explored, but recently zero forcing and power domination for de Bruijn and Kautz digraphs was studied in \cite{GKRS}. De Bruijn and Kautz digraphs are iterated line digraphs of the complete digraph, with and without loops, respectively. In this work, we extend the results in \cite{GKRS} to zero forcing and power domination of iterated line digraphs of any regular digraph.

The line digraph has been used in a broad range of disciplines, but its large number of applications precludes us from including an exhaustive summary here. In the context of this work, since the line digraph of a digraph described by a unitary matrix can also be described by a unitary matrix, iterated line digraphs are used to obtain arbitrarily large digraphs described by unitary matrices {\blu (see \cite{PTZ03} and \cite{S03})}. Such digraphs are frequently used to model quantum systems in physics, chemistry, and engineering {\blu (see \cite{KS97})}, and it was precisely to control quantum systems that zero forcing was introduced in physics. Indeed, line digraphs of digraphs described by unitary matrices are used in quantum computation and in the study of quantum walks {\blu (see \cite{S03})}, as their statistical dynamics models that of random matrix theory {\blu (see \cite{PTZ03} and \cite{S03})}. In particular, the use of regular quantum graphs was studied in \cite{S02} and in \cite{ST04}, as the line digraph of a regular digraph is the digraph of the transition matrix of a coined quantum walk. In addition, line digraphs have been used in information theory as solutions to the {\em index coding with side information problem} {\blu (see \cite{DSC13})}, where the minimum rank of a digraph represents the length of an optimal scalar linear solution of the corresponding instance of the problem {\blu (see \cite{BBJK06})}.

In this work, we extend to digraphs the relationship between zero forcing and power domination established for undirected graphs in \cite{REUF2015}. We also present lower and upper bounds for the zero forcing and the power domination numbers of iterated line digraphs, and show that for regular digraphs, the corresponding lower and upper bounds coincide, providing expressions for the zero forcing and the power domination numbers. Combining our results with known properties of the minimum rank of line digraphs, we conclude that iterated line digraphs of regular digraphs present optimal properties in regards to the minimum rank, zero forcing and power domination problems. We apply our results to the Bruijn and Kautz digraphs, generalized de Bruijn and generalized Kautz digraphs, and wrapped butterflies. Through our work, we show that the relationship between minimum rank, zero forcing and power domination is a powerful tool that permits {\blu us} to combine results obtained separately for each problem to produce stronger results.

\section{Definitions and notation}

A {\em digraph} is a pair $G=(V,A)$, where  $V=V(G)$ is a finite, non-empty set of {\em vertices}, and $A=A(G)$ is a set of ordered pairs of vertices called {\em arcs}. \color{blu}The  {\it order} of $G$ is defined as $|V(G)|$. \color{black} An arc in the form $(u,u)$ is called a {\em loop}. The {\it open out-neighborhood} of a vertex $v$ is $N^+_{G}(v)=\{u\in V(G): (v,u)\in A(G)\}$ and the {\it open in-neighborhood} of $v$ is $N^-_{G}(v)=\{w\in V(G): (w,v)\in A(G)\}$. The {\it closed out-neighborhood} of $v$ is $N^+_{G}[v]=N^+_G(v)\cup \{v\}$ and the {\it closed in-neighborhood} of $v$ is $N^-_{G}[v]=N^-_{G}(v)\cup \{v\}$.  The {\it out-degree}  of $v$ is $d_G^+(v)=|N^+_G(v)|$ and the {\it in-degree} of $v$ is $d_G^-(v)=|N^-_G(v)|$. 
The digraph $G$ is $d${\it-regular} if $d^+_{G}(v)=d^-_{G}(v)=d$, for every $v\in V$. \color{blu}More generally, a digraph $G$ is {\it regular} if there exists an integer $d$ for which $G$ is $d$-regular; otherwise $G$ is said to be {\it irregular}. \color{black}The {\it maximum out-degree} of $G$ is $\Delta ^+(G)=\max \{d_G^+(v):v\in V(G)\}$ and the {\it maximum in-degree} of $G$ is $\Delta ^-(G)=\max \{d_G^-(v):v\in V(G)\}$. The {\it minimum out-degree} of $G$ is  $\delta ^+(G)=\min \{d_G^+(v):v\in V(G)\}$ and the {\it minimum in-degree} of $G$ is  $\delta ^-(G)=\min \{d_G^-(v):v\in V(G)\}$. For a set of vertices $T\subseteq V$,
  $N_G^+[T]=\bigcup_{v\in T}N_G^+[v]$, and analogously for the other neighborhoods. We will omit the subindices when the digraph $G$ is obvious from the context.  
  
  {\color{blu}
If $u$ and $v$ are two different vertices in $G$, a {\em path} of {\it length} $d$ from $u$ to $v$ is a sequence of distinct vertices $u=x_0,\ldots ,x_d=v$ such that $(x_i,x_{i+1})\in A(G)$ for every $i=0,\ldots ,d-1$. A {\em cycle} of {\it length} $\ell$ in $G$, is a sequence of vertices $x_0,\ldots ,x_{\ell}$ such that $x_0,\ldots ,x_{\ell -1}$ are distinct, $x_{\ell}=x_0$ and $(x_i,x_{i+1})\in A(G)$ for every $i=0,\ldots ,\ell-1$.  A digraph $G$ is {\em strongly connected} if for any two vertices $u$ and $v$ there is a path from $u$ to $v$ in $G$.  A digraph is {\em weakly connected} if its underlying graph (i.e. the graph obtained by replacing each arc $(u,v)$ or symmetric pair of arcs $(u,v),(v,u)$ by the edge $uv$) is connected. As suggested by their terms, strong connectivity implies weak connectivity, but they are not equivalent. If a digraph $G$ is not weakly connected, each maximal weakly connected sub-digraph of $G$ is a {\it weak component} of $G$. Note that a digraph with exactly one vertex is weakly connected, so every vertex in $G$ is in exactly one weak component. Thus, the vertex sets of all weak components of $G$ form a partition of the vertex set of $G$.  In a weakly connected digraph that is not strongly connected, the notion of a  {\it strong component} is analogous. For  terminology about graphs or digraphs not defined above we refer the reader to \cite{CLZ}.\\ 

Now we present the notion of zero forcing in digraphs, followed by its formal definition. Intuitively, zero forcing can be described through a coloring process on the vertices of a digraph. Initially, each vertex of a digraph $G$ is arbitrarily colored in one of two colors, say blue and white. Then, apply a given color changing rule that establishes a condition for a white vertex to become blue. Iteratively apply the color changing rule until it fails to produce new blue vertices. At that moment, if all vertices in $G$ are blue, then the initial set of blue vertices is a zero forcing set of $G$. The zero forcing problem consists of finding a zero forcing set of minimum cardinality for a given digraph.  
The color changing rule to be applied in a digraph $G$ depends on whether $G$ has loops or not. For digraphs without loops, the color changing rule is: every white vertex that is the only white out-neighbor of a {\it blue} vertex becomes blue. For digraphs that have at least one loop, the color changing rule is: if a white vertex is the only white out-neighbor of a vertex {\it (blue or white)}, then the white vertex becomes blue. The difference in the color changing rule implies that in a digraph $G$, if there is a loop on a white vertex $v$ and all vertices in $N^+_G(v)\setminus \{v\}$ are blue, then $v$ becomes blue. 

Next, we present the definition of zero forcing using a sequence of sets of vertices  to  describe the blue vertices after each application of the color changing rule. Note that in each application, the color changing rule is simultaneously applied to every white vertex. }

\color{black}
Let $G=(V,A)$ be a digraph. For any \color{blu} non-empty \color{black} set $S\subseteq V(G)$ and any non-negative integer $i$ we define $B^{i} (S)$ by the following rules. 
 \begin{itemize}
 \item[1.]  $B^{0} (S)=S.$
 \item[2.] \color{blu} If $G$ does not have any loops, then for every $i\geq 0$\color{black}\\
 $B^{i+1} (S)=B^{i} (S) \cup \{v  \in  V(G)\setminus B^{i} (S): \exists u\in B^i(S),  N^+_G(u)\setminus B^{i} (S)=\{v\} \}.$
 \item[3.] \color{blu} If $G$ has at least one loop, then for every $i\geq 0$\color{black}\\
 $B^{i+1} (S)=B^{i} (S) \cup \{v  \in  V(G)\setminus B^{i} (S): \exists u\in V,  N^+_G(u)\setminus B^{i} (S)=\{v\} \}.$ 
 \end{itemize}
 
Following \cite{H10}, we say that $S\subseteq V(G)$ is a {\it zero forcing set} of $G$ if there exists a non-negative integer $m$ such that $B^{m} (S)=V(G)$. A {\em minimum zero forcing set} is a zero forcing set of minimum cardinality. The {\em zero forcing number} of $G$ is the cardinality of a minimum zero forcing set and is denoted by $\Z(G)$.  When rule $2$ or rule $3$ is applied, we say that $u$ {\em forces } $v$. Analogously, we say that a set $S\subseteq V$ {\em forces} a set of vertices $W$, when for every $w\in W$ there exists \color{blu} $u\in S$ \color{black} such that $u$ forces $w$. \color{blu} Note that if a digraph $G$ is not weakly connected, then $G$ has $r\geq 2$ weak components $G_1,\ldots , G_r$. In this case, as observed in \cite{ELA13}, $\Z(G)=\sum_{i=1}^r\Z(G_i)$. Since \color{black} zero forcing must be studied independently in each weak component, in this paper we assume all digraphs are at least, weakly connected. \\

For a digraph $G=(V,A)$ of order $n$, the {\em qualitative class} of $G$, i.e. {\blu the matrix family of} $G$,  is the set of matrices $\SG$ defined as  $\SG=\{X\in {\R}^{n\times n} : \mbox{ for } i\not=j, X_{i,j}\not= 0 {\color{blu}\Leftrightarrow}(i,j)\in A(G)\}$ if $G$ does not have loops, and as $\SG=\{X\in {\R}^{n\times n} : X_{i,j}\not= 0 {\color{blu} \Leftrightarrow } (i,j)\in A(G)\}$  if $G$ has at least one loop. {\color{blu}The {\em adjacency matrix} of $G$ is the $n\times n$ matrix $\Ad={\mathscr {A}}(G)$ where  ${\mathscr {A}}_{i,j}=1$ if $(i,j)\in A(G)$ and ${\mathscr {A}}_{i,j}=0$ if $(i,j)\not\in A(G)$. Then,  ${\mathscr {A}}(G)\in\SG$.  }
The {\em maximum nullity} of $G$ is $\M(G) = \max \{\nul X: X \in \sym(G)\}$, and the {\em minimum rank} of $G$ is $\mr(G) = \min \{\rank X : X \in S(G)\}$; clearly $\M(G)+\mr(G)=|V(G)|$. The concept of zero forcing models the process to force zeros in a null vector of a matrix $X\in\SG$, implying  $\M(G)\le\Z(G)$ \cite{H10}. As posed in \cite{AIM08}, it is particularly interesting to identify classes of digraphs $G$ for which  $\M(G)=\Z(G)$. \\

{\color{blu} An important concept in this work is that of power domination. We first present the notion of power domination using a coloring process on the vertices of a digraph, and then give a formal definition. Initially, each vertex of a digraph $G$ is arbitrarily colored either blue or white. In power domination, there are two color changing rules. The first one is applied exactly once at the beginning of the process, and establishes that every white vertex in the out-neighborhood of a blue vertex becomes blue. Then, the coloring process continues with the  application of the second color changing rule of power domination, which coincides with the color changing rule used in zero forcing, and is also iteratively applied until it fails to produce new blue vertices.  At that point, if all vertices in the digraph $G$ are blue, then the original set of blue vertices is a {\it power dominating set} of $G$. The power domination problem consists of finding a power dominating set of minimum cardinality for a given digraph.

Let $G=(V,A)$ be a digraph. For any non-empty set $S\subseteq V(G)$ and any non-negative integer $i$ we define $P^{i} (S)$ by the following rules. 
 \begin{itemize}
 \item[1.]  $P^{0} (S)=S.$
 \item[2.]   $P^{1} (S)=N^+[S].$
 \item[3.] If $G$ does not have any loops, then for every $i\geq 1$\\
 $P^{i+1} (S)=P^{i} (S) \cup \{v  \in  V(G)\setminus P^{i} (S): \exists u\in P^i(S),  N^+_G(u)\setminus P^{i} (S)=\{v\} \}.$
 \item[4.] If $G$ has at least one loop, then for every $i\geq 1$\\
 $P^{i+1} (S)=P^{i} (S) \cup \{v  \in  V(G)\setminus P^{i} (S): \exists u\in V,  N^+_G(u)\setminus P^{i} (S)=\{v\} \}.$ 
 \end{itemize}

We say that $S\subseteq V(G)$ is a {\it power dominating set} of $G$ if there exists a non-negative integer $t$ such that $P^{t} (S)=V(G)$. A {\em minimum power dominating set} is a power dominating set of minimum cardinality. The {\em power domination number} of $G$ is the cardinality of a minimum power dominating set and is denoted by $\gamma_P(G)$.  When rule $2$, $3$ or $4$ is applied, we say that $u$ {\em propagates} to $v$. Analogously, we say that a set $S\subseteq V(G)$ {\em propagates} to a set of vertices $W \subseteq V(G)$, when for every $w\in W$ there exists $u\in S$ such that $u$ propagates to $w$. As in the case of zero forcing, it is sufficient to study power domination on weakly connected digraphs. If a digraph is not weakly connected, then power domination is studied independently in each weak component. 

The power domination problem on digraphs was formally introduced in \cite{A08}, and while not mentioned in the definition itself (\cite[Definition 3.3.1]{A08}),  it is clearly stated in \cite[Pg.4]{A08} that only digraphs without loops were considered. Prior to our work, power domination in digraphs with loops had only been studied in \cite{GKRS}. It is important to remark that in \cite{GKRS}, the authors studied power domination in digraphs with loops using the rules introduced in \cite{A08} for digraphs without loops, while we defined power domination using different rules for digraphs with at least one loop than for digraphs without loops. As a consequence, our definition of power domination differs from the one in \cite{GKRS} in the following situation. If there is a loop on a white vertex $v$ and all vertices in $N^+(v)\setminus\{v\}$ are blue, by our rules $v$ becomes blue, while by the rules in \cite{GKRS} $v$ remains white. By treating digraphs with loops in the same way as in zero forcing, our definition preserves an important relationship between zero forcing and power domination in digraphs without loops, also present in the case of undirected graphs (see \cite{REUF2015}). Indeed, a careful observation of the definition of zero forcing and our definition of power domination in digraphs yields the conclusion that, a set $S\subseteq V(G)$ is a power dominating set of digraph $G$ if and only if $N^+ [S]$ is a zero forcing set of digraph $G$. } \\

For a digraph $G=(V,A)$, the {\it line digraph} of $G$ is the digraph $L(G)$ where $V(L(G))=\{uv: (u,v)\in A(G)\}$ and $A(L(G))=\{uxy:  (u,x), (x,y)\in A(G) \}$. Iterated line digraphs are recursively defined by:  $L^0(G)=G$ and  $L^r(G)=L(L^{r-1}(G))$ for every integer $r\geq 1$. {\color{blu}Observe that in particular, $L^1(G)=L(G)$.} 
Following \cite{P95}, we say that a digraph $G$ is {\em $L$-convergent} if the set $\{L^r:r \mbox{ is a non-negative integer}\}$ is finite; otherwise $G$ is {\em $L$-divergent}. In this paper we are especially interested in digraphs that are $L$-divergent.  In \cite{B68} it  was proven that a digraph $G$ is $L$-divergent if and only if at least one strong component of $G$ is not a cycle or $G$ has at least two cycles joined by a path. The class of $L$-divergent digraphs includes all  strongly connected digraphs other than a cycle, which have been proven to be asymptotically dense in \cite{FYA84} for the regular case and in \cite{D17} for irregular digraphs. \\

In Section 3, we establish lower and upper bounds for the zero forcing number of iterated line digraphs, and determine the zero forcing number of iterated line digraphs of regular digraphs, for which we provide specific constructions of minimum zero forcing sets. We combine the results obtained on zero forcing with known results about minimum rank to prove that iterated line digraphs of regular digraphs are infinite families of digraphs with the property $\M(G)=\Z(G)$. In Section 4, we establish a relationship between power domination and zero forcing in iterated line digraphs and determine the power domination number of iterated line digraphs of regular digraphs. In Section 5, we apply the results to special families of iterated line digraphs.


\section{Zero forcing and minimum rank}

{\color{blu}
We start this section by introducing the definitions of critical and strongly critical sets of vertices in a digraph. Intuitively, both concepts refer to the property of a set of vertices $W$ in a digraph $G$, that if all vertices in $W$ are white and all vertices in $V(G)\setminus W$ are blue, then the color changing rule fails to produce any additional blue vertices. If digraph $G$ does not have loops, this means that no vertex in $V(G)\setminus W$ forces a vertex in $W$, and in this case we say that $W$ is {\it critical}. However, in a digraph with at least one loop, since a vertex could force itself, it is necessary that no vertex, neither in $V(G)\setminus W$  nor  in $W$, can force a vertex in $W$, and in this case, the set $W$ is {\it strongly critical}. 
 }

\begin{defn}
In a digraph $G=(V,A)$, a {\blu non-empty} set {\color{blu}$W\subseteq V(G)$} is called {\em critical} if every $v \in V(G)\setminus W$ has either no {\blu out-neighbors} in $W$, or it has at least {\blu two out-neighbors} in $W$.  That is, for every $v\in V(G)\setminus W$, {\color{blu} $|N^+(v)\cap W|\not=1$}. In addition, $W$ is {\em strongly critical} if for every $v\in V$, {\color{blu}$|N^+(v)\cap W|\not=1$}.
\end{defn}

\begin{rem}\label{atleastone}
Let  $G=(V,A)$ be a digraph and let $S$ be a zero forcing set of $G$. If $G$ does not have any loops, then $|S\cap W|\geq 1$ for every critical set $W$ in $G$.  {\color{blu}Indeed, if $W$ is critical in $G$ and $|S\cap W|=0$, then there is no {\color{blu} $v  \in  V(G)\setminus W$ } such that  $|N^+(v)\cap W|=1$ so no {\color{blu} $v  \in  V(G)\setminus W$ } can force a vertex in $W$. Since $W$ is non-empty and $S$ is a zero forcing set, there must be at least one vertex in $S\cap W$.} Analogously, if $G$ has at least one loop, then $|S\cap W|\geq 1$ for every strongly critical set $W$ in $G$. 
\end{rem}

\begin{obs}\label{more}
Every strongly critical set is a critical set. As a consequence, in any digraph, the maximum number of pairwise disjoint strongly critical sets is less than or equal to the maximum number of pairwise disjoint critical sets. 
\end{obs}

\begin{lem}\label{genlb}
The zero forcing number of a digraph $G$ is at least the maximum number of pairwise disjoint strongly critical sets in $G$. Moreover, if $G$ does not have any loops, its zero forcing number is at least the maximum number of pairwise disjoint critical sets in $G$. 
\end{lem}

\bpf
Let $S$ be a minimum zero forcing set of {\blu digraph} $G$ and let $\{W_1,\ldots ,W_r\}$ be a set of pairwise
disjoint{\blu,} strongly critical sets. By Remark \ref{atleastone}, $|S\cap W_i|\geq 1$ for every $i=1,\ldots , r$, and since the sets $W_1,\ldots ,W_r$ are pairwise disjoint, {\blu hence} $|S|\geq r$. If $G$ does not have any loops, then we apply the same argument with a collection of pairwise disjoint critical sets.
\epf

Next, we will show that the  lower bound provided by Lemma \ref{genlb} can be improved in the case of line digraphs.

\begin{lem}\label{neighbors}
Let $G$ be a digraph and let $uv$ be a vertex of $L(G)$. If  $d^+_{\blu{L(G)}}(uv)\geq 2$, every subset $T\subseteq N^+_{\blu{L(G)}}(uv)$ with $\lvert T\rvert\geq 2$ is a strongly critical set in $L(G)$.
\end{lem}

\bpf
Let $xy$ be any vertex of $L(G)$. If {\blu$y=v$}, then $T\subseteq N^+_{\blu{L(G)}}(xy)$, hence $\lvert N^+_{\blu{L(G)}}(xy)\cap T\rvert=\lvert T\rvert\geq 2$.  If {\blu{$y\not= v$}}, then $|N^+_{\blu{L(G)}}(xy)\cap T|=0$. {\blu Therefore, $T$ is a strongly critical set.}
 \epf

\begin{lem}\label{criticalint}
Let $G$ be a digraph and let $S$ be a zero forcing set of $L(G)$. If $uv$ is a vertex in $L(G)$ such that $d^+_{\blu{L(G)}}(uv)\geq 2$,  then $|N^+_{\blu{L(G)}}(uv)\cap S|\geq d^+_{L(G)}(uv) -1$.

\end{lem}

\bpf
 {\color{blu} 
Let $T=N^+_{L(G)}(uv)\setminus S$, then $|S\cap T|=0$. If $|T|\geq 2$, then $T$ is strongly critical,  by Lemma~\ref{neighbors}, which contradicts Remark~\ref{atleastone}. This implies $|T|\leq 1$, hence $|N^+_{\blu{L(G)}}(uv)\cap S|\geq d^+_{\blu{L(G)}}(uv) -1$.}
\epf

The following simple observation will be used in the proof of the next theorem.

\begin{obs}\label{obs:bad_cycles}
In every digraph $G$, any vertex is in at most one cycle consisting only of vertices with in-degree $1$ in $G$. 
\end{obs}

\begin{thm}\label{linedigraph}
  Let $G=(V,A)$ be a digraph with $\delta ^+(G)\geq 2$ and $\delta ^-(G)\geq 1$. Then, \[ \Z(L(G))= |A(G)|-|V(G)|.\]
\end{thm}

\bpf
Let  $ {\color{blu}V(G)}=\{v_1,\ldots ,v_n\}$ and let $S\subseteq V(L(G))$ be a zero forcing set of $L(G)$. Since $\delta^-(G)\geq 1$, there are vertices $u_1,\ldots,u_n\in V(G)$ such that  $u_iv_i\in  {\color{blu}A(G)}$ for every $i=1,\ldots,n$.  {\color{blu} Then}, by Lemma~\ref{criticalint}, {\color{blu}we have}
  $\lvert S\cap N^+_{\blu{L(G)}}(u_iv_i)\rvert\geq d^+_{\blu{L(G)}}(u_iv_i)-1$ for every $i=1,\ldots,n$. {\color{blu}Further, if $xy\in  N^+_{L(G)}(u_iv_i)\cap N^+_{L(G)}(u_jv_j)$ for some $i\not=j$ then $x=v_i$ and $x=v_j$, a contradiction. Thus, the sets $N^+_{L(G)}(u_iv_i)$ form a partition of $V(L(G))$ meaning $\sum_{i=1}^nd^+_{L(G)}(u_iv_i) = |V(L(G))|=|A(G)|$. Together this gives}
 
  \[\lvert S\rvert=\sum_{i=1}^n\lvert S\cap N^+_{\blu{L(G)}}(u_iv_i)\rvert\geq\sum_{i=1}^n(d^+_{\blu{L(G)}}(u_iv_i)-1)=\lvert
    A(G)\rvert-\lvert V(G)\rvert.\]
  This implies $ {\color{blu} \Z(L(G))}\geq\lvert A(G)\rvert-\lvert V(G)\rvert$.  {\color{blu} To find an upper bound on $\Z(L(G))$ we will now provide a zero forcing set of $L(G)$ of cardinality $|A(G)|-|V(G)|$. }\\
  
 
{\color{blu} From Observation~\ref{obs:bad_cycles}{\blu,} for every $i=1,\ldots ,n$ vertex $v_i$ is in at  most one cycle consisting only of vertices with in-degree $1$, and as a consequence, at most one out-neighbor of $v_i$ is in such cycle. Therefore,  the condition $\delta^+(G)\geq 2$ guarantees that for every $i=1,\ldots,n$ there exists a   vertex $w_i\in N^+_{G}(v_i)$ such that vertex $w_i$ is not in a cycle consisting only of vertices with in-degree $1$ in $G$. If this were not the case, then every $w_i$ would be in such a cycle and arc $v_iw_i\in A(G)$ must also be in that cycle. However, this means $v_i$ is in at least two such cycles, a contradiction. Let}

  \[S=\bigcup_{i=1}^n\{v_iw\ :\ w\in N^+_{\blu{G}}(v_i)\setminus\{w_i\}\}.\] Then, 
  $\lvert S\rvert=\lvert A(G)\rvert-\lvert V(G)\rvert$ and $S$ is a zero forcing set of $G$. Indeed, {\color{blu}for any vertex $v_iv_j$ in $V(L(G))\setminus S$, since  $\delta^-(G)\geq 1$ there exists $v_t\in N^-_G(v_i)$ such that $v_iv_j\in N^+_{L(G)}(v_tv_i)$. If $v_tv_i\in S$, then $v_iv_j=v_iw_i$ so $v_iv_j$ is the only out-neighbor of $v_tv_i$ not forced, by construction of $S$, and $v_tv_i$ forces $v_iv_j$. If $v_tv_i\not\in S$, then proceed with $v_tv_i$ as we did with $v_iv_j$, and the selection of vertices $w_1,\ldots ,w_n$ guarantees that at some point, a vertex in $S$ is reached.}
\epf

\begin{cor}\label{bounds}
Let $G=(V,A)$ be a digraph with $\delta ^+(G)\geq 2$ and $\delta ^-(G)\geq 1$.  Then,  \[ ( \delta  -1 ){\color{blu}|V(G)|}\leq \Z(L(G))\leq   ( \Delta  -1 ){\color{blu}|V(G)|}\]
  
\noindent where $\delta = \max \{ \delta ^+(G), \delta ^-(G) \} $ and $\Delta = \min \{ \Delta ^+ (G), \Delta^-(G)\}. $
\end{cor}

\bpf
From Theorem \ref{linedigraph}, $\Z(L(G))={\color{blu} |A(G)|-|V(G)|}$. Since ${\color{blu}|A(G)|}=\sum _{v\in V}d^+_G(v)$,  then $ \delta ^+(G){\color{blu}|V(G)|}\leq {\color{blu}|A(G)|}\leq \Delta ^+(G){\color{blu}|V(G)|}$ and  $( \delta ^+(G)-1){\color{blu}|V(G)|}\leq \Z(L(G))\leq ({\color{blu}\Delta ^+(G) }-1){\color{blu}|V(G)|}$. Analogously, ${\color{blu}|A(G)|}=\sum _{v\in V}d^-_G(v)$ implies $ (\delta ^-(G)-1){\color{blu}|V(G)|}\leq \Z(L(G))\leq (\Delta ^-(G)-1){\color{blu}|V(G)|}$, and the result is obtained by choosing the {\blu least} upper bound and the {\blu greatest }lower bound.
\epf

The next results extend Theorem \ref{linedigraph} and Corollary \ref{bounds} to iterated line digraphs. 

\begin{thm}\label{itlinedigraph}
Let $G=(V,A)$ be a digraph with {\color{blu}$\delta ^+(G)\geq 2$ and $\delta ^-(G)\geq 1$.} For every integer $n\geq 1$, $ \Z(L^{n}(G))= |V(L^{n}(G))|-|V(L^{n-1}(G))|$.
\end{thm}

\bpf
From Theorem \ref{linedigraph}, $\Z(L^{n}(G))= |A(L^{n-1}(G))|-|V(L^{n-1}(G))|$. By definition of line digraph, $|A(L^{n-1}(G))|=|V(L^{n}(G))|$ so  $\Z(L^{n}(G))= |V(L^{n}(G))|-|V(L^{n-1}(G))|$. 
\epf

{\color{blu}
\begin{lem}\label{degit1}
Let $G=(V,A)$ be a digraph with $\delta ^+(G)\geq 1$ and $\delta ^-(G)\geq 1$. For any integer $n\geq 1$, \\
   \noi \mbox{  }\hskip 1.5cm  a) $\delta^+(L^n(G))=\delta^+(G)$. \hskip 3cm b) $\Delta^+(L^n(G))=\Delta^+(G)$.\\
   \noi \mbox{  }\hskip 1.5cm  c) $\delta^-(L^n(G))=\delta^-(G)$.  \hskip 3cm d) $\Delta^-(L^n(G))=\Delta^-(G)$.
\end{lem}

\bpf
Observe that since $L^n(G)=L(L^{n-1}(G))$, it is sufficient to prove the statements for $n=1$. By definition of line digraph, if $uv\in V(L(G))$, then $d^+_{L(G)}(uv)=d^+_G(v)$ and $d^-_{L(G)}(uv)=d^-_G(u)$. Therefore,  $\delta^+(L(G))\geq \delta^+(G)$, 
 $\delta^-(L(G))\geq \delta^-(G)$,  $\Delta^+(L(G))\leq \Delta^+(G)$ and $\Delta^-(L(G))\leq \Delta^-(G)$. If  $\delta ^-(G)\geq 1$, then for every  $v\in V(G)$ there exists $u\in N^-_G(v)$, and as a consequence, there exists $uv\in V(L(G))$ such that $d^+_{L(G)}(uv)=d^+_G(v)$. Therefore, $ \delta^+(G)\geq \delta^+(L(G))$, $\Delta^+(L(G))\geq \Delta^+(G)$ and we obtain a) and b).  Analogously, $\delta ^+(G)\geq 1$ implies that for every $u\in V(G)$ there exists $uv\in V(L(G))$ such that $d^-_{L(G)}(uv)=d^-_G(u)$. Thus, $ \delta^-(G)\geq \delta^-(L(G))$,  $\Delta^-(L(G))\geq \Delta^-(G)$ and we obtain c) and d). 
\epf

\begin{lem} \label{degit2}
Let $G=(V,A)$ be a digraph with $\delta ^+(G)\geq 1$ and $\delta ^-(G)\geq 1$. For any integer $n\geq 0$,  
\begin{itemize}
\item[a)] $\delta^+(G) |V(L^{n}(G))| \leq |V(L^{n+1}(G))| \leq |V(L^{n}(G))| \Delta^+(G)$.
\item[b)] $\delta^-(G)|V(L^{n}(G))| \leq  |V(L^{n+1}(G))| \leq |V(L^{n}(G))| \Delta^-(G).$
\end{itemize}
\end{lem}

\bpf
Since $ |V(L^{n+1}(G))|= |A(L^{n}(G))|$ and $|A(L^{n}(G))|=\sum_{v\in V(L^{n}(G))} d^+_{L^{n}(G)}(v)$, from the observation 
 \[\delta^+(L^{n}(G)) |V(L^{n}(G))|\mbox{ } \leq \sum_{v\in V(L^{n}(G))} d^+_{L^{n}(G)}(v) \mbox{ }\leq \mbox{ }|V(L^{n}(G))| \Delta^+(L^{n}(G)),\]
\noi we conclude $\delta^+(L^{n}(G)) |V(L^{n}(G))|\leq |V(L^{n+1}(G))| \leq |V(L^{n}(G))| \Delta^+(L^{n}(G)).$ 

By Lemma \ref{degit1}, $\delta^+(L^{n}(G))=\delta^+(G)$ and $\delta^-(L^{n}(G))=\delta^-(G)$, and as a consequence, the previous inequality can be written as $\delta^+(G) |V(L^{n}(G))|\leq  |V(L^{n+1}(G))| \leq |V(L^{n}(G))| \Delta^+(G)$. 

Analogously, from $ |V(L^{n+1}(G))|= |A(L^{n}(G))|$ and $|A(L^{n}(G))|=\sum_{v\in V(L^{n}(G))} d^-_{L^{n}(G)}(v)$, we conclude
 $\delta^-(L^{n}(G)) |V(L^{n}(G))|\leq |V(L^{n+1}(G))| \leq |V(L^{n}(G))| \Delta^-(L^{n}(G))$.  Since Lemma \ref{degit1} also implies $\delta^-(L^{n}(G))=\delta^-(G)$ and $\Delta^-(L^{n}(G))=\Delta^-(G)$, the previous inequality is equivalent to $\delta^-(G)|V(L^{n}(G))| \leq  |V(L^{n+1}(G))| \leq |V(L^{n}(G))| \Delta^-(G)$.
\epf

\begin{lem}\label{degit}
Let $G=(V,A)$ be a digraph with $\delta ^+(G)\geq 1$ and $\delta ^-(G)\geq 1$. For any  integer $n\geq 1$, \begin{itemize}
\item[a)] $(\delta^+(G))^{n} |V(G)| \leq |V(L^{n}(G))| \leq |V(G)| (\Delta^+(G))^{n}$.
\item[b)] $(\delta^-(G))^{n} |V(G)|\leq  |V(L^{n}(G))| \leq |V(G)| (\Delta^-(G))^{n}.$
\end{itemize}
\end{lem}

\bpf 
For $n=1$ the result holds, since $\delta^+(G)|V(G)| \leq |V(L(G))| \leq |V(G)| \Delta^+(G)$ and $\delta^-(G)|V(G)| \leq |V(L(G))| \leq |V(G)| \Delta^-(G)$ follow from replacing $n=0$ in Lemma \ref{degit2}. We conclude the proof by induction on $n$.

 Assume $(\delta^+(G))^{n} |V(G)| \leq |V(L^{n}(G))| \leq |V(G)| (\Delta^+(G))^{n}$ for an integer $n\geq 1$. Then,

\noi $(\delta^+(G))^{n+1} |V(G)| \leq \delta^+(G)|V(L^{n}(G))| \leq |V(L^{n}(G))|\Delta^+(G)\leq |V(G)| (\Delta^+(G))^{n+1}.$ \mbox{ }\hfill \eqnum\label{multiply}\\

\noi By Lemma \ref{degit2}, $\delta^+(G) |V(L^{n}(G))| \leq |V(L^{n+1}(G))| \leq |V(L^{n}(G))| \Delta^+(G),$  \mbox{ }\hfill \eqnum\label{lemma}\\

\noi Combining \eqref{multiply} and \eqref{lemma}, we obtain $(\delta^+(G))^{n+1} |V(G)| \leq |V(L^{n+1}(G))| \leq |V(G)| (\Delta^+(G))^{n+1}$.\\

Now, assume $(\delta^-(G))^{n} |V(G)| \leq |V(L^{n}(G))| \leq |V(G)| (\Delta^-(G))^{n}$ for an integer $n\geq 1$. 
\noi Observe that Lemma \ref{degit2} implies $\delta^-(G) |V(L^{n}(G))| \leq |V(L^{n+1}(G))| \leq |V(L^{n}(G))| \Delta^-(G)$. Therefore, repeating the previous 
argument we conclude $(\delta^-(G))^{n+1} |V(G)| \leq |V(L^{n+1}(G))| \leq  |V(G)| (\Delta^-(G))^{n+1}$.
\epf
}

\begin{cor}\label{itbounds}
Let $G=(V,A)$ be a digraph with {\color{blu} $\delta ^+(G)\geq 2$ and $\delta ^-(G)\geq 1$. }For every integer $n\geq 1$, \[ (\delta -1)\delta ^{n-1}{\color{blu}|V(G)|}\leq \Z(L^n(G))\leq (\Delta -1)\Delta ^{n-1}{\color{blu}|V(G)|},\] 
\noindent where $\delta = \max \{ \delta ^+(G), \delta ^-(G) \} $ and $\Delta = \min \{ \Delta ^+ (G), \Delta^-(G)\}. $
\end{cor}

\bpf
{\blu Lemma \ref{degit2} gives two lower bounds and two upper bounds for $|V(L^{n}(G))|$. Choosing the greatest lower bound and the least upper bound we obtain $\delta |V(L^{n-1}(G))| \leq |V(L^{n}(G))| \leq |V(L^{n-1}(G))| \Delta$. This chain of inequalities implies $(\delta -1)|V(L^{n-1}(G))| \leq |V(L^{n}(G))| -|V(L^{n-1}(G))|\leq |V(L^{n-1}(G))| (\Delta -1)$, and since by Theorem \ref{itlinedigraph}, $\Z(L^n(G))= |V(L^{n}(G))|-|V(L^{n-1}(G))|$, we conclude
\begin{center}
\noi \mbox{  }\hskip 1.5cm $(\delta -1)|V(L^{n-1}(G))| \leq \Z(L^n(G))\leq |V(L^{n-1}(G))| (\Delta -1).$ \mbox{ }\hfill \eqnum\label{zfton-1}
\end{center}

By Lemma \ref{degit}, again, selecting the greatest lower bound and the least upper bound, we have
 $\delta^{n-1} |V(G)|\leq |V(L^{n-1}(G))|\leq  |V(G)|\Delta^{n-1}$. From these inequalities we obtain
$(\delta-1)\delta ^{n-1}|V(G)|\leq(\delta-1)|V(L^{n-1}(G))|$ and $ |V(L^{n-1}(G))|(\Delta-1)\leq |V(G)|\Delta ^{n-1}(\Delta-1),$ which combined with \eqref{zfton-1} yield $(\delta -1)\delta ^{n-1}|V(G)|\leq \Z(L^n(G))\leq (\Delta -1)\Delta ^{n-1}|V(G)|.$}
\epf

In the {\blu remainder} of this section we restrict ourselves to regular digraphs. The next result follows immediately from Theorem \ref{itlinedigraph} or from Corollary \ref{itbounds}.

\begin{cor}\label{itregular}
Let  $G=(V,A)$ be a $d$-regular digraph with $d\geq 2$.  For every integer $n\geq 1$, 
\begin{itemize}{\color{blu}
\item [a)]  $\Z(L^{n}(G))=(d-1)d^{n-1}{\color{blu}|V(G)|}.$
\item [b)]  $\Z(L^{n+1}(G))=d\Z(L^{n}(G)).$}
\end{itemize}
\end{cor}

The following result, together with our results on zero forcing of line digraphs, allow us to show that for any $d$-regular digraph $G$ and any positive integer $n$,   $\M(L^n(G))=\Z(L^n(G))$.

\begin{lem}\cite[Lemma 1]{GW02}
\label{Gimbert}
A $d$-regular digraph of order $p$ is a line digraph if and only if the rank of its adjacency matrix is equal to $p/d$.
\end{lem}

\begin{thm}
\label{thm1}
Let $G=(V,A)$ be a $d$-regular digraph with $d\geq 2$. For every integer $n\geq 1$, 
\begin{itemize}{\color{blu}
\item [a)]  $\M(L^{n}(G))=\Z(L^{n}(G))=(d-1)d^{n-1}{\color{blu}|V(G)|}.$
\item [b)]  $\mr (L^{n}(G))=d^{n-1}{\color{blu}|V(G)|}.$}
\end{itemize}
\end{thm}

\bpf
Let $\Ad=\Ad(L^n(G))$ {\blu be} the adjacency matrix of $L^n(G)$.  Since $G$ is $d$-regular, then $L^n(G)$ is also $d$-regular, and by Lemma \ref{Gimbert}, $\Ad$ has rank ${{d^n|V(G)|}\over {d}}=d^{n-1}|V(G)|$. Then, the equality  ${\color{blu} \nul(\Ad)+\rank(\Ad)}=|V(L^n(G))|=d^n|V(G)|$ implies that $\nul(\Ad)=d^n|V(G)|-d^{n-1}|V(G)|=(d-1)d^{n-1}|V(G)|.$ Since $\nul(\Ad)\leq \M(L^{n}(G))\leq \Z(L^{n}(G))$ and, by Corollary \ref{itregular}, $\Z(L^{n}(G))=(d-1)d^{n-1}|V(G)|$, we conclude $(d-1)d^{n-1}|V(G)|=\nul(\Ad)\leq \M(L^{n}(G))\leq \Z(L^{n}(G))=(d-1)d^{n-1}|V(G)|$ and the result follows immediately. Replacing this value in $\mr (L^{n}(G))=d^n|V(G)|-\M(L^{n}(G))$ completes the proof.
\epf

\begin{rem}
The only $1$-regular{\blu, weakly connected,} digraphs are cycles. If $G$ is a cycle, then $L^n(G)=G$ for any integer {\blu$n\geq 1$}, so {\color{blu} $\Z(G)=1$, $\mr (G)=|V(G)|-1$, $\M(G)=1$.}  As a consequence, for any regular digraph $G$ and for any integer {\blu$n\geq 1$},   $\M(L^{n}(G))=\Z(L^{n}(G))$. Thus, iterated line digraphs of regular digraphs have optimal values of zero forcing, minimum rank and maximum nullity. Furthermore, the minimum rank and maximum nullity are attained by the adjacency matrix of $L^{n}(G)$.
\end{rem}

\section{Power domination}


Let {\blu$G$} be a digraph and {\blu$S$ a set of vertices of $G$}. {\blu{As observed in Section 2,}} $S$ is a power dominating set of $G$ if and only if $N^+_{\blu{G}}[S]$ is a zero forcing set of $G$. Hence, $\Z(G)\leq \gamma_P(G)(\Delta ^+(G)+1)$. The following result provides an improved upper bound in the case that $G$ is a line digraph. 
\begin{thm}\label{pdlb}
Let $G$ be a digraph with $\delta ^+(G)\geq 1$. Then, \[ \Z(L(G))\leq \Delta ^+(G)\gamma_P(L(G)). \] 
\end{thm}


\bpf
Let  $S=\{u_1v_1,\ldots ,u_rv_r\}$ be a minimum power dominating set of $L(G)$. We prove  $\Z(L(G))\leq \Delta ^+(G)\gamma_P(L(G))$ by constructing a zero forcing set of $L(G)$ with cardinality $\Delta ^+(G)|S|$. 

{\blu{D}}efine $S_1=\{u_iv_i\in S\,:\,v_i\neq v_j\text{ for all }j\text{ with }1\leq j<i\}$ and $S_2=S\setminus S_1$. Since $S$ is a power dominating set of $L(G)$,   $N^+_{\blu{L(G)}}[S]=\cup_{i=1}^r N^+_{\blu{L(G)}}[u_iv_i]=S_2 \cup{\blu{\{}}\cup _{uv\in S_1} N^+_{\blu{L(G)}}[uv] {\blu{\}}}$ is a zero forcing set of $L(G)$.  Since $\delta^+(G)\geq 2$, for every $u_iv_i\in S_1$ we can select an arbitrary vertex $v_iw_i\in N^+_{\blu{L(G)}}(u_iv_i)$. Let $P=\cup_{u_iv_i\in S_1} (N^+_{\blu{L(G)}}[u_iv_i] \setminus \{v_iw_i\})$.  Then{\blu,} $P$ forces $\cup _{uv\in S_1} N^+_{\blu{L(G)}}[uv]$, so $S_2\cup P$ forces $S_2 \cup{\blu{\{}}\cup _{uv\in S_1}N^+_{\blu{L(G)}}[uv]{\blu{\}}}=N^+_{\blu{L(G)}}[S]$. Since $N^+_{\blu{L(G)}}[S]$ is a zero forcing set of $L(G)$, then $S_2\cup P$ is also a zero forcing set of $L(G)$ and $\Z(G)\leq |S_2\cup P|$. Now, $|S_2\cup P|\leq  |S_2|+\sum  _{u_iv_i\in S_1} |N^+_{\blu{L(G)}}[u_iv_i]\setminus \{w_iu_i\})|\leq  |S_2|+|S_1|\Delta^+(G)\leq  |S|\Delta^+(G)$. Therefore, $\Z(L(G))\leq \gamma_P(L(G))\Delta ^+(G)$.
\epf


\begin{cor}\label{supergood}
Let $G$ be a digraph with {\color{blu} $\delta ^+(G)\geq 2$ and $\delta ^-(G)\geq 1$}.  For every integer $n\geq 1$, 
 \[ \Z(L^n(G))\geq \gamma_P(L^{n}(G)) \geq \blu{\left\lceil {{\Z(L^{n}(G))}\over {\Delta ^+(G)}}\right\rceil.} \] 
\end{cor}

\begin{cor}\label{regularbound}
For a regular digraph $G$, $\Z(L(G))\leq \gamma_P(L^{2}(G))$.
\end{cor}

\bpf
 Assume $G$ is $d$-regular.  Then{\blu,} by definition of line digraph, $L(G)$ is also $d$-regular. If $d=1$, then {\color{blu}$G=L^n(G)$ are cycles,} so $\Z(L(G))=1$ and $\gamma_P(L^{2}(G))=1$. If $d\geq 2$, we apply Theorem \ref{pdlb} to $L(G)$ and obtain  $\Z(L^2(G))\leq d \gamma_P(L^2(G))$. By Corollary \ref{itregular}{\blu,} $\Z(L^2(G))=d\Z(L(G))$ so $d\Z(L(G)) \leq d \gamma_P(L^2(G))$, and this implies $\Z(L(G)) \leq  \gamma_P(L^2(G))$.
\epf


Next, we will show that if $G$ is a $d$-regular digraph, then $\Z(L(G))= \gamma_P(L^2(G))$. First, we need to introduce additional terminology and obtain further results.

{\color{blu}
\begin{dfn}
A digraph $H$ is a {\it factor} of a digraph $G$ if $V(H)=V(G)$ and $A(H)\subseteq A(G)$. In particular, if $H$ is a factor of $G$ and $H$ is $1$-regular, then $H$ is a $1${\it-factor} of $G$. 
\end{dfn}
\begin{obs}
A $1$-factor in a digraph $G$ is either a cycle, or a set of pairwise vertex-disjoint cycles, containing every vertex in $G$. 
\end{obs} }
\begin{thm}\label{good}
Let $G$ be a digraph with {\color{blu} $\delta ^+(G)\geq 2$ and $\delta^-(G)\geq 1$.  If $G$ has a $1$-factor in which every cycle contains at least one vertex $v$ with  $d^-_G(v)>1$,} then 
$\Z(L(G))\geq \gamma_P(L^{2}(G))$. 
\end{thm}
\bpf
{\color{blu} Let $H$ be a $1$-factor of $G$ with the condition in the hypothesis.} Then, $H$ induces a permutation $f$ on the set of vertices of $G$ where $f(v)=u$ if and only if $(u,v)$ is an arc in $H$. Observe that the cycles in $H$ correspond to the orbits of the permutation $f$.  Since each cycle in $H$ contains a vertex  $v$ with  $d^-_G(v)>1$, then each orbit in the permutation $f$ contains a vertex  $v$ with  $d^-_G(v)>1$. By definition, $V(L(G))=\{uv: (u,v)\in A(G)\}$ and $V(L^2(G))=\{uvw: (u,v)\in A(G)\mbox{ and } (v,w)\in A(G)\}$. Define $S=\{ f(u)uv \in V(L^2(G)): u\not= f(v)\}$. Observe  $|S|=\sum_{v\in V(H)} (d^-_{G}(v)-1)$ and $V(H)=V(G)$, therefore  $|S|=\sum_{v\in V(G)} (d^-_G(v)-1).$
Besides, $\sum_{v\in V(G)} (d^-_G(v)-1) =(\sum_{v\in V(G)} d^-_G(v))-|V(G)| = |A(G)|-|V(G)|$. By Theorem \ref{linedigraph}, $\Z(L(G))=|A(G)|-|V(G)|$, and as a consequence, it is sufficient to prove that $S$ is a power dominating set of $L^2(G)$ to conclude $\gamma_P(L^{2}(G))\leq \Z(L(G))$.\\

Given $xyz$ in $V(L^2(G))$, it is sufficient to show that either $xyz\in N^+_{L^2(G)}[S]$ or $xyz$ is obtained by a sequence of forces starting with a vertex in $N^+_{L^2(G)}[S]$.
  \begin{itemize} 
  \item[ ] {\it Case 1.} If $x\neq f(y)$ then $f(x)xy\in S$, hence $xyz\in N^+_{L^2(G)}(f(x)xy)\subseteq N^+_{L^2(G)}[S]$.
  \item[ ] {\it Case 2.} If $x=f(y)$ and $y\neq f(z)$ then $xyz\in S\subseteq N^+_{L^2(G)}[S]$.
  \item[ ] {\it Case 3.} If $x=f(y)$ and $y= f(z)$ then $xyz=f^2(z)f(z)z$. 
  \end{itemize}
  
 In the last case,  for notational convenience, denote the elements of the orbit of $z$ by $z_i=f^{i}(z)$, i.e., $z=z_0=z_\ell$,
    $y=z_1=z_{\ell+1}$, etc., where $\ell$ is the length of the orbit of $z$. By assumption,
    $d^-(z_k)\geq 2$ for some $k\geq 2$. This implies that there exists $a\in N^-_G(z_k)$ such that
    $a\neq z_{k+1}$. Consequently, $f(a)az_k\in S$ and
    $az_kz_{k-1}\in N^+_{L^2(G)}(f(a)az_k)\subseteq N^+_{L^2(G)}[S]$. If $d^+(z_{k-1})=1$ then $az_kz_{k-1}$ forces its only out-neighbor $z_kz_{k-1}z_{k-2}$. If $d^+_G(z_{k-1})\geq 2$ every out-neighbor of $az_kz_{k-1}$ different from
        $z_kz_{k-1}z_{k-2}$ is in the form $z_kz_{k-1}b$ for some vertex
        $b\in N^+_G(z_{k-1})\setminus \{z_{k-2}\}$. Since $f$ is bijective and $z_{k-2}\neq b$, we
        have $z_{k-1}\neq f(b)$, and therefore, $z_kz_{k-1}b \in S$. So every out-neighbor of
        $az_kz_{k-1}$ different from $z_kz_{k-1}z_{k-2}$ is in $S$ and this means that $az_kz_{k-1}$
        forces $z_kz_{k-1}z_{k-2}$. Repeating this argument, we obtain that all the vertices $z_{i}z_{i-1}z_{i-2}$ for
      $i=k,k-1,\dots ,2$ are forced, and in particular $xyz=z_2z_1z_0$.
\epf
\color{blu}
  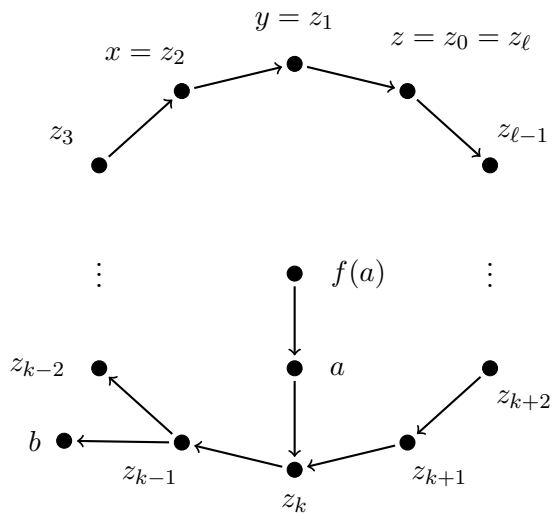
\begin{figure}[htb]
    \centering
    \begin{tikzpicture}[yscale=.9,every node/.style={draw,circle,fill=black,outer sep=2pt,inner
        sep=2pt,minimum size=0pt}]
      \node[label={[label distance=-.3cm]60:{\small $z=z_0=z_\ell$}}] (z0) at (60:3) {};
      \node[label={[label distance=-.3cm]90:{\small $y=z_1$}}] (z1) at (90:3) {};
      \node[label={[label distance=-.2cm]120:{\small $x=z_2$}}] (z2) at (120:3) {};
      \node[label={[label distance=.1cm]150:{\small $z_3$}}] (z3) at (150:3) {};
      \node[label={[label distance=-.1cm]30:{\small $z_{\ell-1}$}}] (zl-1) at (30:3) {};
      \node[label={[label distance=-.1cm]330:{\small $z_{k+2}$}}] (zk+2) at (330:3) {};
      \node[label={[label distance=-.1cm]300:{\small $z_{k+1}$}}] (zk+1) at (300:3) {};
      \node[label={[label distance=-.1cm]270:{\small $z_k$}}] (zk) at (270:3) {};
      \node[label={[label distance=-.1cm]240:{\small $z_{k-1}$}}] (zk-1) at (240:3) {};
      \node[label={[label distance=.1cm]180:{\small $z_{k-2}$}}] (zk-2) at (210:3) {};
      \node[label={[label distance=.1cm]0:{\small $a$}}] (a) at (270:1.5) {};
      \node[label={[label distance=.1cm]0:{\small $f(a)$}}] (fa) at (270:.1) {};
      \node[label={[label distance=-.1cm]180:{\small $b$}}] (b) at (220:4) {};
      \draw[thick,->] (z3) -- (z2);
      \draw[thick,->] (z2) -- (z1);
      \draw[thick,->] (z1) -- (z0);
      \draw[thick,->] (z0) -- (zl-1);
      \draw[thick,->] (zk+2) -- (zk+1);
      \draw[thick,->] (zk+1) -- (zk);
      \draw[thick,->] (zk) -- (zk-1);
      \draw[thick,->] (zk-1) -- (zk-2);
      \draw[thick,->] (fa) -- (a);
      \draw[thick,->] (a) -- (zk);
      \draw[thick,->] (zk-1) -- (b);
      \path (zk-2) -- node[draw=none,fill=none]{\vdots} (z3);
      \path (zl-1) -- node[draw=none,fill=none]{\vdots} (zk+2);      
    \end{tikzpicture}    
    \caption{{\blu Illustration of the argument for Case 3 in the proof of Theorem 4.7.}}
    \label{fig:proof_thm_4.4}
  \end{figure}
\color{black}


{\blu In a digraph $G$ with $\delta^-(G)\geq 2$, every $1$-factor satisfies the condition in Theorem \ref{good}. As a consequence, we obtain the following corollary.}

\begin{cor} \label{supergood2}
Let $G$ be a digraph with $\delta ^+(G)\geq 2$ and {\color{blu}$\delta^-(G)\geq 2$.} If $G$ has a $1$-factor, then $Z(L(G))\geq \gamma_P(L^{2}(G)).$
\end{cor}
{\color{blu}



Next, we recall another definition and a property of regular digraphs from \cite{S06}.
\begin{dfn}
A {\em cycle factorization} in a digraph $G$ is a set of $1$-factors $\{H_1,\ldots, H_k\}$ such that $A(G)=\cup^k_{i=1} A(H_i)$ and $A(H_i)\cap A(H_j)=\emptyset$ for every $1\leq i<j\leq k$.  
 \end{dfn}
 
\begin{lem}\cite[Lemma 1]{S06}\label{dfact}
Let $G$ be a $d$-regular digraph. Then, $G$ has a cycle factorization consisting of exactly $d$ $1$-factors.
\end{lem}

By Lemma \ref{dfact}, every regular digraph has a cycle factorization, and as a consequence, it has at a $1$-factor. Therefore, when $G$ is $d$-regular and $d\geq 2$ combining Corollary \ref{regularbound} and Corollary \ref{supergood2} we obtain the following result.

\begin{thm}\label{partial1}
If $G$ is a regular digraph, then $Z(L(G))= \gamma_P(L^{2}(G)).$
\end{thm}

\bpf
Assume $G$ is $d$-regular. If $d=1$, then as seen in the proof of Corollary \ref{regularbound}, $\Z(L(G))=1$ and $\gamma_P(L^{2}(G))=1$. If $d\geq 2$, then $\Z(L(G))\leq \gamma_P(L^{2}(G))$, by Corollary \ref{regularbound}. By Lemma \ref{dfact}, $G$ has a cycle factorization. As as a consequence, $G$ has a $1$-factor, and by Corollary \ref{supergood2}, $\Z(L(G))\geq \gamma_P(L^{2}(G))$. Thus, $\Z(L(G))= \gamma_P(L^{2}(G)).$
\epf




We now extend Theorem \ref{partial1} to a relationship between $\Z(L^n(G))$ and $\gamma_P(L^{n+1}(G))$ for any regular digraph $G$ and any positive integer $n$.
}
\begin{cor}\label{secondld2}
Let $G$ be a $d$-regular digraph. For any integer $n\geq 1$, 
\[ \Z(L^n(G))= \gamma_P(L^{n+1}(G)).\]
Moreover, if $d\geq 2$, then {\blu$\Z(L^n(G))= \gamma_P(L^{n+1}(G))=(d-1)d^{n-1}|V(G)|$.}
\end{cor}

\bpf
{\color{blu}If $d=1$, then $G$ is a cycle so $G=L^n(G)$ for every $n\geq 1$ and  it follows that $\Z(L^n(G))= \gamma_P(L^{n+1}(G))=1$. Assume $d\geq 2$. Observe that since $G$ is $d$-regular, then  $L^n(G)$ is also $d$-regular, for every integer $n\geq 1$.  By Corollary \ref{partial1}, $\Z(L^n(G))= \gamma_P(L^{n+1}(G))$ holds for $n=1$.  For $n\geq 2$, we apply Corollary \ref{partial1} to $L^{n-1}(G)$ and obtain  $\Z(L(L^{n-1}(G)))= \gamma_P(L^{2}(L^{n-1}(G)))$. Observe that 
$L(L^{n-1}(G))=L^n(G)$ and $L^{2}(L^{n-1}(G))=L(L^n(G))=L^{n+1}(G)$. Then, $\Z(L(L^{n-1}(G)))= \gamma_P(L^{2}(L^{n-1}(G)))$ implies $Z(L^n(G))= \gamma_P(L^{n+1}(G))$.   By Corollary \ref{itregular}, if $d\geq 2$, then $ \Z(L^n(G))= {(d-1)d^{n-1}|V(L(G))|}$.  Since we have proven $Z(L^n(G))= \gamma_P(L^{n+1}(G))$, we conclude that $ \gamma_P(L^{n+1}(G))= (d-1)d^{n-1}|V(G)|.$  }\epf

\smallskip

{\color{blu}
Observe that Theorem \ref{partial1} does not provide an expression to determine $\gamma_P(L(G))$. We show next that $\gamma_P(L(G))$ depends on properties of the digraph $G$, and provide an expression to determine $\gamma_P(L(G))$ if $G$ is a $d$-regular digraph. Since the only $1$-regular digraphs are cycles, we assume $d\geq 2$, and when regularity is not necessary, $\delta ^+(G)\geq 2$ and $\delta ^-(G)\geq 2$. 

\begin{lem}\label{unique}
Let $G$ be a digraph with $\delta ^+(G)\geq 2$ and $\delta ^-(G)\geq 2$. If there exists a set $S\subseteq V(G)$ such that for any two different vertices $x,y\in S$, $N^+_G(x)\cap N^+_G(y)=\emptyset$,
then, for each $v\in N^+_G(S)$ there exists a unique vertex $u\in S\cap N^-_G(v)$.
\end{lem}

\bpf
Suppose $v\in N^+_G(S)\setminus S$ and there exist two different vertices $u,w\in S$ such that  $u\in N^-_G(v)$ and $v\in N^-_G(v)$. Then,  $v\in N^+_G(u)\cap N^+_G(w)$, a contradiction. Now suppose $v\in N^+_G(S)\cap S$ and consider two cases depending on $(v,v)\in A(G)$ or $(v,v)\not\in A(G)$. If $(v,v)\in A(G)$, then $v\in N^-_G(v)$. Suppose there exists $u\not=v$ such that $u\in N^-_G(v)$. Then, $N^+_G(u)\cap N^+_G(u)=\{v\}$, a contradicting. If $(v,v)\not\in A(G)$ and there exist two different vertices $u,w\in S$ such that  $u\in N^-_G(v)$ and $v\in N^-_G(v)$. Then,  $v\in N^+_G(u)\cap N^+_G(w)$, another contradiction.
\epf

\begin{prop}\label{ubisn-gen}
Let $G$ be a digraph with $\delta ^+(G)\geq 2$ and $\delta ^-(G)\geq 2$. If there exists a set $S\subseteq V(G)$ such that for any two different vertices $x,y\in S$, $N^+_G(x)\cap N^+_G(y)=\emptyset$, and for every $x\in S$, $N^+_G(x)\cap S=\emptyset$ or $N^+_G(x)\cap S=\{x\}$, then  $\gamma_P(L(G))\leq |V(G)|-|S|$. 
\end{prop}

\bpf
Assume $V(G)=\{v_1,\ldots ,v_n\}$ and $S=\{v_r,\ldots, v_n\}$. For each integer $i=1,\ldots ,n$ select one vertex $u_i\in N^-_G(v_i)$ following this rule: if $v_i\in N^+_G(S)$, then select as $u_i$ the unique vertex, by Lemma \ref{unique}, $u_i\in (S\cap N^-_G(v_i))$. Define $P=\{u_iv_i: 1\leq i\leq r-1\}.$ Observe that $V(L(G))=\cup_{i=1}^n N^+_{L(G)}(u_iv_i)$ and $N^+_{L(G)} [P]=P\cup (\cup_{i=1}^{r-1}N^+_{L(G)}(u_iv_i))$. 

If $v_i\in S$ and there is no loop $(v_i,v_i)\in A(G)$, then the selection of vertices $u_i$ implies that for each vertex $v_j\in N^+_G(v_i)$, the vertex selected as $u_j$ is $v_i$. As a consequence, $v_jv_i\in P$ and $N^+_{L(G)}(u_iv_i)\subseteq P$. If there is a loop  $(v_i,v_i)\in A(G)$, then $v_iv_i\in V(L(G))$ and $v_iv_i\not\in P$. However, by the rule for selecting $u_i$, $(N^+_{L(G)}(v_iv_i)\setminus \{v_iv_i\})\subset P$. Then,  $v_iv_i$ forces itself, and we conclude that $P$ forces all vertices in $N^+_{L(G)}(v_iv_i)$. 
Then $N^+_{L(G)}[P]$ is a zero forcing set of $L(G)$, $P$ is a power dominating set of $L(G)$, and since $|P|=|V(G)|-|S|$, we conclude $\gamma_P(L(G))\leq |V(G)|-|S|$. 
\epf

In any digraph $G$ with $\delta ^+(G)\geq 2$ and $\delta ^-(G)\geq 2$, for any $v\in V(G)$ with $d^+_G(v)$, the set $S=\{v\}$ satisfies the conditions of Proposition \ref{ubisn-gen} so $\gamma_P(L(G))= |V(G)|-1$. 

\begin{cor}\label{ubis1}
Let $G$ be a $d$-regular digraph of order $n$. If $d\geq 2$ and $n<2d$, then $\gamma_P(L(G))= |V(G)|-1. $
\end{cor}

\bpf
By Theorem \ref{pdlb}, $\left\lceil {{\Z(L(G))}\over d} \right\rceil \leq \gamma_P(L(G))$. By Corollary \ref{itregular}, $\Z(L(G))=n(d-1)$.  Then, $n-\left\lfloor {{n}\over d} \right\rfloor \leq \gamma_P(L(G))$. Therefore, if $\left\lfloor {{n}\over d} \right\rfloor=1$, we conclude $\gamma_P(L(G))= |V(G)|-1$. 
\epf

\begin{cor}\label{ubisn-gen}
Let $G$ be a $d$-regular digraph of order $n$. If $d\geq 2$ and there exists $S\subseteq V(G)$ such that 1) for any two different vertices $x,y\in S$, $N^+_G(x)\cap N^+_G(y)=\emptyset$, 2) for every $x\in S$, $N^+_G(x)\cap S=\emptyset$ or $N^+_G(x)\cap S=\{x\}$, and 3) $|S|= \left\lfloor {n\over d}\right\rfloor$, then  $\gamma_P(L(G))= \left\lceil {{n(d-1)}\over d}\right\rceil$. 
\end{cor}

}

\section{Applications}\label{Applications}

Next, we recall the definitions of some families of iterated line digraphs that have been
extensively used in applications. In each family, we apply Corollary \ref{itregular} to obtain their
zero forcing number, Theorem \ref{thm1} to obtain their maximum nullity and their minimum rank.{\blu
  The power domination number of $L^n(G)$ follows from Corollary \ref{secondld2} when $n\geq 2$, and
  from Corollary \ref{ubis1} or Corollary \ref{ubisn-gen} when $n=1$.} We refer the reader to
\cite{G} for additional details on the families of digraphs studied in this section.
 
 \subsection{de Bruijn digraphs}
 
 For any integers $d\geq 2$ and $D\geq 1$, the {\em de Bruijn digraph} $B(d,D)$ has for vertices the set $\mathbb {Z}_{d^D}$, and each vertex $x$ is adjacent to the vertices $dx+t$ for any $ t\in \mathbb {Z}_d$. Alternatively,  $B(d,D)$ can be iteratively defined by following rules: $B(d,1)=K_d$, the complete symmetric digraph of order $d$ with a loop on each vertex, and $B(d,D)=L(B(d,D-1))$, if $D\geq 2$. That is, $B(d,D)=L^{D-1}(K_d)$, so $B(d,D)$ has order $d^D$.

\begin{cor}
For any integers $d\geq  2$ {\blu{and}} $D\geq 2$:

\begin{enumerate}
\item $\Z(B(d,D))=(d-1)d^{D-1}$. 
\item $\M(B(d,D))=(d-1)d^{D-1}$.
\item $\mr(B(d,D))=d^{D-1}$.
\item $\gamma_P(B(d,D))=(d-1)d^{D-2}$.
\end{enumerate}

\end{cor}

{\blu
For $D=2$, $B(d,2)=L(K_d)$ and $\gamma_P(B(d,2))=d-1$ follows from Corollary \ref{ubis1}.}\\

Quantum systems based on de Bruijn digraphs were studied in \cite{PTZ03} and \cite{S03}.  Random walks on de Bruijn digraphs $B(d,n)$ have the fastest mixing rates among $d$-regular digraphs of order $d^n$ \cite{S03}.

\subsection{Kautz digraphs} 
 
For any integers $d\geq 3$ and $D\geq 1$, the vertices of the {\em Kautz digraph} $K(d,D)$ are the $D$-tuples{\blu in $\mathbb {Z}_{d+1}$} in which any two consecutive elements are different; each vertex $x$ is adjacent to the vertices $dx+t$ for any $ t\in \mathbb {Z}_d$. The digraphs $K(d,D)$ can also be iteratively defined by following rules: $K(d,1)=K^*_{d+1}$,  the complete symmetric digraph of order $d+1$ without loops, and $K(d,D)=L(K(d,D-1))$, if $D\geq 2$. Therefore,  $K(d,D)=L^{D-1}(K^*_{d+1})$,  and $K(d,D)$ has order $d^{D-1}(d+1)$.

\begin{cor}
For any integers $d\geq 3$ and $D\geq 3$: 

\begin{enumerate}
\item $\Z(K(d,D))=(d-1)d^{D-2}{(d+1)}$.
\item $\M(K(d,D))=(d-1)d^{D-2}{(d+1)}$.
\item $\mr(K(d,D))=d^{D-2}{(d+1)}$.
\item $\gamma_P(K(d,D))=(d-1)d^{D-3}{(d+1)}$.
\end{enumerate}
\end{cor}

 \subsection{Generalized de Bruijn digraphs} 

For any integers $d\geq 2$ and $n\geq 1$, the {\em generalized de Bruijn digraph} $GB(d,n)$, has $\mathbb {Z}_n$ as its vertex set, and each vertex $x$ is adjacent to vertices $(dx+t)\mod n$ for any value of $t$ in $\mathbb {Z}_d$. The name of this family of digraphs is due to the fact that it {\blu contains} the de Bruijn digraphs as a sub-family. More precisely, $B(d,D)=GB(d,d^D)$. The generalized de Bruijn digraphs, also known as {\em Reedy-Pradhan-Kuhl digraphs}, have the property that $GB(d,dn)=L(GB(d,n))$, so for every integer $m\geq 1$, $GB(d,d^mn)=L^m(GB(d,n))$.

\begin{cor}
For any integers $d\geq 2$, $n\geq 1$ and {\blu$m\geq 2$:}

\begin{enumerate}
\item $\Z(GB(d,d^mn))=(d-1)d^{m-1}n$.
\item $\M(GB(d,d^mn))=(d-1)d^{m-1}n$.
\item $\mr(GB(d,d^mn))=d^{m-1}n$.
\item $\gamma_P(GB(d,d^mn))=(d-1)d^{m-2}n$.
\end{enumerate}
\end{cor}

  \subsection{Generalized Kautz digraphs} 

For any integers $d\geq 2$ and $n\geq 1$, the  {\em generalized Kautz digraph} $GK(d,n)$ has $\mathbb {Z}_n$ as its vertex set, and each vertex $x$ is adjacent to vertices $(-dx-t)\mod n$ for any value of $t$ in $\mathbb {Z}_d$. This family contains the Kautz digraphs as a sub-family. Indeed, $K(d,D)=GK(d,d^D+d^{D-1})$. The generalized Kautz digraphs, also known as {\em Imase-Itoh digraphs}, have the following property with respect to the line digraph: $GK(d,dn)=L(GK(d,n))$. Therefore,  for every integer $m\geq 1$, $GK(d,d^mn)=L^m(GK(d,n))$. 

\begin{cor}
For any integers $d\geq 2$, $n\geq 1$ and {\blu$m\geq 2$:}

\begin{enumerate}
\item $\Z(GK(d,d^mn))=(d-1){\color{blu}d^{m-1}}n$.
\item $\M(GK(d,d^mn))=(d-1){\color{blu}d^{m-1}}n$.
\item $\mr(GK(d,d^mn))=d^{m-1}n$.
\item $\gamma_P(Z(GK(d,d^mn)))=(d-1)d^{m-2}n$.
\end{enumerate}
\end{cor}

 \subsection{Directed wrapped butterfly} 

For any  integers $d,n\geq 2$, the {\em directed wrapped butterfly} $WB(d,n)$ has for vertices the ordered pairs $(x,l)$ where $x$ is an  $n$-tuple of integers in $\mathbb{Z}_d$ and $l$ is an integer, $0\leq l\leq n-1$. A vertex $(x_1\ldots x_n, l)$ is adjacent to $d$ vertices in the form  $(x_1\ldots x_{l-1}\alpha x_{l+1} x_n, l)$ where $\alpha \in \mathbb{Z}_d$. The wrapped butterflies $WB(d,n)$ can be obtained as a line digraph, but we need some definitions prior to explain this in detail. 

If $G$ and $H$ are two digraphs, the {\em conjunction} of $G$ and $H$ is the digraph $G\otimes H$, whose vertex set corresponds to $V(G)\times V(H)$; a vertex $(g,h)$ is adjacent to a vertex $(g',h')$ if $(g,g')$ is an arc in $G$ and $(h,h')$ is an arc in $H$. Then, $WB(d,n)=L^{n-1}(K_d\otimes C_n)$ where $K_d$ denotes the complete symmetric digraph of order $d$ with a loop on each vertex, and $C_n$ denotes the cycle of order $n$. 

\begin{cor}
For any integers $d\geq 2$ and $n\geq 2$:

\begin{enumerate}
\item $\Z(WB(d,n))=(d-1)d^{n-1}n$.
\item $\M(WB(d,n))=(d-1)d^{n-1}n$.
\item $\mr(WB(d,{\color{blu}n}))=d^{n-1}n$.
\item $\gamma_P(WB(d,n))=(d-1)d^{n-2}n$.
\end{enumerate}
\end{cor}

{\blu
For $n=2$, $WB(d,2)=(K_d\otimes C_2)$ and $\gamma_P(WB(d,2))=2(d-1)$ follows from Corollary\ref{ubisn-gen}.}\\

\subsection*{Acknowledgment} The authors are grateful for the thoughtful referee comments that helped them to substantially improve the content and  presentation of this article.


\begin{thebibliography}{20}

\bibitem{A08} A. Aazami, {\em Hardness results and approximation algorithms for some problems on graphs,} Ph.D. thesis, University of Waterloo, 2008.

\bibitem{AIM08} AIM Minimum Rank -- Special Graphs Work Group (F. Barioli, W. Barrett, S. Butler, S. M. Cioab\u{a}, D. Cvetkovi\'c, S. M. Fallat, C. Godsil, W. Haemers, L. Hogben,  R. Mikkelson,  S. Narayan,  O. Pryporova,   I. Sciriha,  W. So,   D. Stevanovi\'c,  H. van der Holst, K. Vander Meulen,  A. Wangsness),  {\em Zero forcing sets and the minimum rank  of graphs,}   Linear Algebra App. 428 (2008) 1628--1648.

\bibitem{G} J. Bang-Jensen, G.Z. Gutin, {\em Digraphs - Theory, Algorithms and Applications}, Springer-Verlag, London, 2009.

\bibitem{BBJK06} Z. Bar-Yossef, Z. Birk, T. S. Jayram, T. Kol, {\em Index coding with side information,} IEEE Trans. Inform. Theory. 57, (2011) 1479--1494.

\bibitem{B09} F. Barioli, S. Fallat, D. Hershkowitz, H.T. Hall, L. Hogben, H. van der Holst, B. Shader,  {\em On the minimum rank of not necessarily symmetric matrices: a preliminary study,} Electron. J. Linear Algebra. 18 (2009), 126--{\blu145.}

\bibitem{B68} L.W. Beineke, {\em On derived graphs and digraphs}, in H. Sachs (Ed.), Beitraege zur Graphentheorie, Teubner-Verlag, Leipzig (1968), 17--23.

\bibitem{REUF2015} 	K.F. Benson, D. Ferrero, M. Flagg, V. Furst, L. Hogben, V. Vasilevska, B. Wissman, {\em Zero forcing and power domination for graph products,} {\blu  Australas. J. Combin. 70 (2018) 221--235.  }

\bibitem{ELA13} A. Berliner, M. Catral, L. Hogben, M. Huynh, K. Lied, M. Young,{\em  Minimum rank, maximum nullity and zero forcing number of simple digraphs,} Electron.
J. Linear Algebra 26 (2013) 762--780.

\bibitem{BG07} D. Burgarth, V. Giovannetti,  {\em Full control by locally induced relaxation,}  Phys. Rev. Lett.  99 (2007) 100501.

\bibitem{BFS99} J.F. Buss,  G.S. Frandsen,  J.O. Shallit,  {\em The  computational complexity of some problems of linear algebra.}  J. Comput. System Sci. 58 (1999) 572--596.

\bibitem{CLZ}{\blu G. Chartrand, L. Lesniak, P. Zhang, {\em Graphs \& Digraphs}, 6th. ed.,  CRC Press, Boca Raton, 2016.}

\bibitem{D17} C. Dalf\'o, {\em Iterated line digraphs are asymptotically dense,} Linear Algebra Appl. 529 (2017) 391--396.

\bibitem{DSC13} S.H. Dau, V. Skachek, Y.M.  Chee,  {\em Error correction for index coding with side information,}   IEEE Trans. Inform. Theory. 59 (2013) 1517--1531.

\bibitem{DR} P. A. Dreyer Jr., Fred S. Roberts, {\em Irreversible $k$-threshold processes: Graph-theoretic threshold models of the spread of disease and of opinion,} Discrete Appl. Math. 157 (2009) 1615--{\blu 1627.}

\bibitem{FH11} S. Fallat, L. Hogben, {\em Variants on the minimum rank problem: A survey II,}  arXiv:102.5142, 2011.

\bibitem{FHB04} M. Fazel, H. Hindi, S. Boyd,  {\em Rank minimization and applications in system theory, } Proceedings of the American Control Conference (2004) 3273--3278.

\bibitem{FYA84} M.A. Fiol, J.L.A. Yebra, I. Alegre, {\em Line digraphs iterations and the $(d,k)$ digraph problem,} IEEE Trans. Comput. 33 (1984) 400--403.

\bibitem{GW02} J. Gimbert, Y. Wu, {\em The underlying line digraph structure of some $(0,1)$-matrix equations,} Discrete Appl. Math. 116 (2002) 289--296. 

\bibitem{GKRS} C. Grigorious, T. Kalinowski, S. Stephen, {\em On the power domination number of de Bruijn and Kautz digraphs,} {\blu  In: Brankovic L., Ryan J., Smyth W. (eds.) Combinatorial Algorithms IWOCA 2017, LNCS 10765, Springer, 2018.}


\bibitem{HHHH02}  T.W. Haynes, S.M. Hedetniemi, S.T. Hedetniemi, M.A. Henning, {\em Domination in graphs applied to electric power networks,} SIAM J. Discrete Math. 15 (2002) 519--529.

\bibitem{H10} L. Hogben,  {\em Minimum rank problems,}  Linear Algebra Appl. 432 (2010)1961--1974.

\bibitem{KS97} T. Kottos, U. Smilansky, {\em Quantum chaos on graphs,}  Phys. Rev. Lett. 79 (1997) 4794--4797.

\bibitem{MZC14} N. Monshizadeh, S. Zhang, M. K. Camlibel, {\em  Zero forcing sets and controllability of dynamical systems defined on graphs,} IEEE Trans. Automat. Control 59 (2014) 2562--2567.

\bibitem{PTZ03} P. Pak\'onski, G. Tanner, K. \.{Z}yczkowski, {\em Families of line-graphs and their quantization,} J. Stat. Phys 111 (2003) 1331--1352.

\bibitem{P76} T. D. Parsons, {\em Pursuit-evasion in a graph,} Theory and Applications of Graphs. Springer-Verlag, 426--441, 1976.

\bibitem{P95} E. Prisner, {\em Graph dynamics,} Longman House, Harlow (UK), 1995. 

\bibitem{S06} S. Severini, {\em On the structure of the adjacency matrix of the line digraph of a regular digraph,} Discrete Appl. Math. 154 (2006) 1663--1665.

\bibitem{S03} S. Severini, {\em  On the digraph of a unitary matrix,} SIAM J. Matrix Anal. Appl. 25 (2003) 295--300.

\bibitem{S02} S. Severini, {\em The underlying digraph of a coined quantum random walk,} arXiv:0210055, 2002.

\bibitem{ST04} S. Severini, G. Tanner, {\em Regular quantum graphs,}  J. Phys. A: Math. Gen. 37 (2004) 6675--6686.

\bibitem{TD15} M. Trefois, J-C. Delvenne, {\em Zero forcing number, constrained matchings and strong structural controllability,}   Linear Algebra Appl. 484 (2015) 199--218.

\bibitem{B13} B. Yang, {\em Fast--mixed searching and related problems on graphs,} Theoret. Comput. Sci. 507 (2013) 100--113.

\end{thebibliography}
\end{document}